\documentclass[12pt]{amsart}
\usepackage{amsfonts,amsmath,amsthm,amssymb,}
\usepackage[usenames,dvipsnames,svgnames,table]{xcolor}
\usepackage{nicematrix}
\usepackage{tikz}

\usepackage{fullpage} 

\usepackage{todonotes}

\numberwithin{equation}{section}

\newtheorem{thm}{Theorem}[section]
\newtheorem*{thm*}{Theorem}

\newtheorem*{thm*A}{Theorem A}

\newtheorem{prop}[thm]{Proposition}
\newtheorem*{prop*}{Proposition}

\newtheorem*{question*}{Question}

\newtheorem{lemma}[thm]{Lemma}

\newtheorem*{conj*}{Conjecture}
\theoremstyle{definition}
\newtheorem{defin}[thm]{Definition}

\newtheorem{example}[thm]{Example}
\newtheorem{remark}[thm]{Remark}
\newtheorem*{remark*}{Remark}
\newtheorem*{remarks*}{Remarks}

\newcommand{\ip}[1]{\langle #1 \rangle}

\newcommand{\myspan}{\operatorname{span}}
\newcommand{\Rspan}{\mathbb R\text{-}\myspan}
\newcommand{\Aut}{\operatorname{Aut}}
\newcommand{\mytrace}{\operatorname{trace}}

\newcommand{\Der}{\operatorname{Der}}

\newcommand{\Id}{\operatorname{Id}}

\begin{document}

\title{Concerning  a conjecture of Taketomi-Tamaru}
\author[Jablonski]{Michael Jablonski}
\thanks{This work was supported in part by NSF grant DMS-1906351.}
\date{December 31, 2024}
\maketitle

\begin{abstract}  We study the setting of 2-step nilpotent Lie groups in the particular case that its type $(p,q)$ is not exceptional.  We demonstrate that, generically, the orbits of $\mathbb R^{>0}\times \Aut_0$ in $GL(n)/O(n)$ are congruent even when a Ricci soliton metric does exists.  In doing so, we provide a counterexample to the local version of a conjecture of Taketomi-Tamaru.  
\end{abstract}

Among solvable and nilpotent groups, perhaps the most natural distinguished Riemannian metrics are those left-invariant metrics which are either Einstein or Ricci soliton.  These metrics are known to minimize natural functionals \cite{Heber,Lauret:RicciSolitonHomogeneousNilmanifolds}, have maximal symmetry when compared to other left-invariant metrics \cite{Jablo:ConceringExistenceOfEinstein,
Jablo:MaximalSymmetryAndUnimodularSolvmanifolds,
GordonJablonski:EinsteinSolvmanifoldsHaveMaximalSymmetry,
GordonJablonski:RicciSolitonSolvmanifoldsHaveInfinitesimalMaximalSymmetry}
, and constitute the full class of non-compact, homogeneous Einstein and Ricci soliton metrics \cite{Bohm-Lafuente:HomogeneousEinsteinMetricsOnEuclideanSpacesAreEinsteinSolvmanifolds}.
  The pursuit of algebraic and geometric criteria which guarantee or preclude the existence of these metrics has been a long standing avenue of investigation, with the lion's share of the attention given to the algebraic side, see e.g.\ \cite{Nikolayevsky:EinsteinSolvmanifoldsandPreEinsteinDerivation}.

Recently there have been several works approaching this question from the geometric side.  Given a Lie group $G$ with Lie algebra $\mathfrak g$, one can study the left-invariant metrics on $G$ by studying inner products on $\mathfrak g$; thus, one considers the set of inner products which is naturally presented as the symmetric space $GL(n)/O(n)$, where $n=\dim \mathfrak g$.

The subgroup $\mathbb R^*\times \Aut(\mathfrak g) \subset GL(n)$ acts on the space of inner products $GL(n)/O(n)$.  Any two inner products in the same $\Aut(\mathfrak g)$-orbit are isometric and so inner products in the same $\mathbb R^*\times \Aut(\mathfrak g)$-orbit are isometric up to scaling.  Note, in the nilpotent setting the $\Aut(\mathfrak g)$-orbits are precisely the isometry classes, but this is not necessarily true for solvable, non-nilpotent groups \cite{GordonWilson:IsomGrpsOfRiemSolv}.

Given an inner product $\ip{\cdot, \cdot } \in GL(n)/O(n)$, the orbit $\mathbb R^*\times \Aut(\mathfrak g) \cdot \ip{\cdot, \cdot }$ has recently been dubbed the \textit{corresponding submanifold}.  We collect some recent results on the geometry of this corresponding submanifold.

\begin{itemize}
	\item  Let $G$ be a 3-dimensional solvable Lie group with left-invariant metric $\ip{\cdot, \cdot } \in GL(3)/O(3)$.  Then $(G,\ip{\cdot, \cdot })$ is a Ricci soliton if and only if the corresponding submanifold $\mathbb R^*\times \Aut(\mathfrak g)\cdot \ip{\cdot, \cdot }$ in $GL(3)/O(3)$ is a minimal submanifold 
\cite{Hashinaga-Tamaru:Three-dimensionalSolvsolitonsAndTheMinimalityOfTheCorrespondingSubmanifolds}.

	\item Let $G$ be a 4-dimensional nilpotent Lie group with left-invariant metric $\ip{\cdot, \cdot }\in GL(4)/O(4)$.  Then $(G,\ip{\cdot, \cdot })$ is a Ricci soliton if and only if the corresponding submanifold $\mathbb R^*\times \Aut(\mathfrak g)\cdot \ip{\cdot, \cdot }$ in  $GL(4)/O(4)$ is a minimal submanifold.  Furthermore, it is shown that this result is false for some solvable groups in dimension 4 
\cite{Hashinaga-OnTheMinimalityOfTheCorrespondingSubmanifoldsToFour-DimensionalSolvsolitons}.

	\item Let $G$ be an n-dimensional solvable Lie group.  In \cite{Taketomi:OnARiemannianSubmanifoldWhoseSliceRepresentation}, a sufficient condition for a metric to be Ricci soliton is given in terms of the so-called slice representation of $\mathbb R^*\times \Aut(\mathfrak g)$ acting on $GL(n)/O(n)$.

	\item Let $G$ be an n-dimensional solvable Lie group.  It is conjectured that if $\mathbb R^*\times \Aut(\mathfrak g)$ does not act transitively on $GL(n)/O(n)$ and all the orbits are congruent, then $G$ does not admit a Ricci soliton \cite[Conj.~1.2]{Taketomi-Tamaru:OnTheExistenceOfLeftInvariantRicciSolitons-AConjectureAndExamples}.  The authors go on to verify the conjecture for some classes of nilpotent groups in every dimension - the algebras there have large automorphism groups.

\end{itemize}

Recall, two orbits being congruent means that there is an isometry $\phi \in GL(n)$ of the symmetric space $GL(n)/O(n)$ which sends one orbit to the other.

In the present work, we  consider the connected components of the corresponding submanifolds, i.e. we will look at the orbits of the connected group $\mathbb R^{>0}\times \Aut(\mathfrak g)_0$ in $GL(n)/O(n)$.

\begin{thm*A}\label{thm: main theorem}
There exists a 9-dimensional nilpotent Lie group $G$ such that
	\begin{enumerate}
	\item $G$ admits a Ricci soliton metric and
	\item $\mathbb R^{>0}\times \Aut(\mathfrak g)_0$ does not act transitively on  $GL(9)/O(9)$ and the orbits are all congruent.
	\end{enumerate}
\end{thm*A}
See Example \ref{example: 9-dimensional nilpotent} for details.  This result provides a counterexample to the local version of the conjecture of Taketomi-Tamaru - i.e., on the connected components of the corresponding submanifolds.   Whether or not the full conjecture holds remains an open and interesting question.  

Even further, our result demonstrates  that  there  is no  criterion on the local geometry of the corresponding submanifolds which can completely determine the soliton condition, or likely any other distinguished condition.  Whether or not the soliton condition can be characterized by the global action of $\mathbb R^*\times \Aut(\mathfrak g)$  on $GL(n)/O(n)$ is an open and interesting question.  For partial progress on this question, see  \cite{Taketomi:OnARiemannianSubmanifoldWhoseSliceRepresentation}.

We note that we do not take up the remaining interesting question of whether or not  it is possible that conditions such as the corresponding submanifold  being minimal or not either guarantee or preclude the existence of a soliton metric.  It would be interesting to know if the  corresponding submanifolds appearing in Theorem A are indeed minimal submanifolds.

Finally, we point out that the phenomenon occurring in the theorem above is not special.  In fact, this is the behavior seen in the generic setting, cf.\  Theorem \ref{thm: generic points being soliton and small derivations}.

\section{outline of proof} \label{sec: outline of proof}\label{sec: outline of proof}

Our strategy is to work with Lie algebras whose algebra of derivations is very small.  Recall, the Lie algebra of the automorphism group is the algebra of derivations.  

Consider a 2-step nilpotent Lie algebra $\mathfrak n = \mathfrak v + \mathfrak z$, where $\mathfrak z=[\mathfrak n,\mathfrak n]$ is the commutator subalgebra and $\mathfrak v$ is a complement of $\mathfrak z$.  As $\mathfrak n$ is 2-step nilpotent, $\mathfrak z=[\mathfrak n,\mathfrak n]$ is central, but might not be all of the center of $\mathfrak n$.  The algebra $\mathfrak n$ comes equipped naturally with two kinds of derivations.  First we have the (1,2)-derivation 
	\begin{equation}\label{eqn: 1,2 derivation}
		 D =  \begin{bmatrix}
	\Id_\mathfrak v & 0\\
	0 & 2 \ \Id_\mathfrak z \end{bmatrix}	   ;
	\end{equation} 
then we have derivations which vanish on $\mathfrak z$ and map $\mathfrak v $ to $\mathfrak z$, i.e. ones of the form
	\begin{equation} \label{eqn: der v to z}
	 \begin{bmatrix}
	0 & 0\\
	* & 0   \end{bmatrix}	.
	\end{equation}
These are precisely the derivations valued in $\mathfrak z$ and they form an ideal in $\Der(\mathfrak n)$.  We denote this set of derivations by $\Der_{\mathfrak v \to \mathfrak z}$.  Combining  the above, we have a subalgebra $\mathbb R (D) \ltimes \Der_{\mathfrak v \to\mathfrak z}$ of $\Der(\mathfrak n)$.  In fact, one can argue that this is an ideal of $\Der(\mathfrak n)$, though we won't need this fact.    In general, the set of derivations is
	\begin{equation} \label{eqn: derivations split}
	\Der(\mathfrak n) = \Der(\mathfrak n)\cap (\mathfrak{gl(v)}\oplus \mathfrak{gl(z)}) \oplus \Der_{\mathfrak v \to \mathfrak z};
	\end{equation}
and at the automorphism level one has
	$$\Aut(\mathfrak n) = (\Aut(\mathfrak n) \cap GL(\mathfrak v)\times GL(\mathfrak z)) \ \exp (\Der_{\mathfrak v\to\mathfrak z}),$$
where the subgroup $\exp(\Der_{\mathfrak v\to\mathfrak z})$ is a normal subgroup.

We are interested in the case when $\Der(\mathfrak n)$ is as small as possible; i.e.,  $\Der(\mathfrak n) = \mathbb R (D) \ltimes \Der_{\mathfrak v \to\mathfrak z}$.  For the rest of this section we make this assumption going forward.  As we will see, this is the generic situation, cf.~Theorem \ref{thm: generic derivations}.

By Lie's Theorem, we know there is some upper triangular matrix algebra $\mathfrak s \supset \Der(\mathfrak n)$.  When  $\Der(\mathfrak n)$ is minimal, this is easily seen using a basis of $\mathfrak v$ concatenated with a basis of $\mathfrak z$.  Observe,  when  $\Der(\mathfrak n)$ is minimal, it is then an ideal of $\mathfrak s$.  If $S$ denotes the subgroup of $GL(n,\mathbb R)$ with Lie algebra $\mathfrak s$ then we have
	\begin{itemize}
	\item $\mathbb R^{>0}\times \Aut(\mathfrak n)_0$ is a normal, proper subgroup of $S$ and
	\item $S$ acts transitively on $GL(n,\mathbb R)/O(n)$.
	\end{itemize}
Using  normality of $\mathbb R^{>0}\times\Aut(\mathfrak n)_0$ and transitivity of $S$ above, we immediately see that all the $\mathbb R^{>0}\times\Aut(\mathfrak n)_0$-orbits in $GL(n,\mathbb R)/O(n)$ are congruent; furthermore, as  $\mathbb R^{>0}\times \Aut(\mathfrak n)_0$ is a proper subgroup of $S$, it does not act transitively on $GL(n)/O(n)$.

Below we detail for which types of 2-step nilpotent Lie groups the above arguments come together.  In a generic sense, many algebras simultaneously 
	\begin{itemize}
	\item 	admit a soliton metric and 
	\item have smallest possible derivation algebra.
\end{itemize}		
These two facts combined give the existence of a counterexample.  However, putting your finger on a counterexample is significantly more challenging than one might example - see Example \ref{example: 9-dimensional nilpotent}.

\section{The $j$-map and algebras of type $(p,q)$}
Consider a 2-step nilpotent algebra $\mathfrak n = \mathfrak v + \mathfrak z$ where $\mathfrak z = [\mathfrak n , \mathfrak n]$.  We say $\mathfrak n$ is of type $(p,q)$ if $\dim \mathfrak z = p$ and $\dim \mathfrak v = q$.  Given an algebra of type $(p,q)$, with an inner product $\ip{\cdot, \cdot }$, we may consider the so-called $j$-map studied by Eberlein and others \cite{Eberlein:2step}:
	$$j:\mathfrak z \to \mathfrak{so(v)},$$
defined by 
	$$\ip{j(z)v,w} = \ip{[v,w],z}.$$
Taking an orthonormal basis $\{z_1,\dots, z_p\}$ of $\mathfrak z$, one may associate to $\mathfrak n$ a $p$-tuple of skew-symmetric matrices
	$$C = (C_1, \dots, C_p) \in \mathfrak{so}(q) ^p = \mathfrak{so}(q) \otimes \mathbb R^p$$
via $C_i = j(z_i) \in \mathfrak{so(v)}$.  The set of 2-step nilpotent algebras of type $(p,q)$ forms a Zariski open set $V_{pq}^0$ in $\mathfrak{so}(q) \otimes \mathbb R^p$, being those $p$-tuples whose entries are linearly independent.  We note that the constraint of linear independence forces us to have $1\leq p \leq \frac 1 2 q (q-1) = \dim \mathfrak{so}(q)$.

Interestingly, both the automorphism group and the ismorphism classes of 2-step nilpotent algebras of type $(p,q)$ can be read off from a natural $GL(q)\times GL(p)$ action on $\mathfrak{so}(q) \otimes \mathbb R^p$.  This action is given as follows.  For $(g,h)\in GL(q)\times GL(p)$ and $M\otimes v\in\mathfrak{so}(q)\otimes \mathbb R^p$,
	$$(g,h)\cdot M\otimes v = gMg^t\otimes hv,$$
where the $GL(p)$ action on $\mathbb R^p$ is the standard one; of course, one extends linearly.  We note that there is an induced Lie algebra action of $\mathfrak{gl}(q)\oplus \mathfrak{gl}(p)$ given by
	\begin{equation}\label{eqn: gl x gl Lie algebra action}
	(X,Y) \cdot M\otimes v = (XM + MX^t)\otimes v + M\otimes Yv,
	\end{equation}
for $X\in\mathfrak{gl}(q)$ and $Y\in\mathfrak{gl}(p)$.  This action and its relationship to nilpotent geometry have been explored in depth by Eberlein \cite{EberleinModuli}.  We record some useful facts here, cf. Equation \ref{eqn: derivations split}.
	\begin{itemize}
	\item Two algebras of type $(p,q)$ are isomorphic if and only if their corresponding $p$-tuples lie in the same $GL(q)\times GL(p)$ orbit.
	\item If $\mathfrak n$ corresponds to $C\in\mathfrak{so}(q)\otimes \mathbb R^p$, then 
		$$\Der(\mathfrak n)\cap (\mathfrak{gl}(q)\oplus \mathfrak{gl}(p)) \simeq  (\mathfrak{gl}(q)\oplus \mathfrak{gl}(p))_C  ,$$
		where the right-hand side is the stabilizer of the Lie algebra action at $C$.
	\end{itemize}
	
The isomorphism $\Psi : \Der(\mathfrak n)\cap (\mathfrak{gl}(q)\oplus \mathfrak{gl}(p)) \to  (\mathfrak{gl}(q)\oplus \mathfrak{gl}(p))_C$ above is given by
	\begin{align}\label{eqn: isomorphism between derivations and stabilizer subalgebra}
	\Psi (X,Y) = (-X^t,Y).
	\end{align}
At the group level, we have an isomorphism $\Aut(\mathfrak n)\cap (GL(q)\times GL(p)) \to  (GL(q)\times GL(p))_C$ given by $\Psi (g,h) = \left( \left(g^t\right)^{-1}, h\right)$.

\begin{defin}\label{definition of exception pairs (p,q)}
We say that $(p,q)$ is exceptional if either $(p,q)$ or $(\frac 1 2 q (q-1) - p, q)$ appears in the list below.  Note, $\frac 1 2 q(q-1) = \dim \mathfrak{so}(q)$.
	\begin{quote}
	$(1,q)$ for $q \geq 2$\\ 
	$(\frac 1 2 q(q-1), q)$ for $q \geq 2$\\ 
	$(2,k)$ for $k \geq 3$\\
	$(3,k)$ for $4\leq k\leq 6$
	\end{quote}
\end{defin}

As we will see below, the exceptional types $(p,q)$ are when the generic derivation algebras are larger than the minimal possible one appearing in the next theorem.  Interestingly, most of the go-to examples that people work with fall into the exceptional cases above and so do not reflect the nature of generic 2-step nilpotent geometry.

\begin{thm}  \label{thm: generic derivations}
For non-exceptional types, a generic algebra has as its derivation algebra the minimal possible derivation algebra, i.e.
	$$\Der =  \mathbb R(D) \oplus \Der_{\mathfrak v \to \mathfrak z},$$
where $D$ is the $(1,2)$-derivation given in Eqn.~\ref{eqn: 1,2 derivation}.
More precisely, there exists a Zariski open set in $\mathfrak{so}(q)^p$ with the property above.
\end{thm}

\begin{lemma}  Take $C\in \mathfrak{so}(q)^p$.  If $SL(q)\times SL(p) \cdot C$ is closed, then the stabilizer of the $\mathfrak{gl}(q)\oplus \mathfrak{gl}(p)$ action at $C$ is $\mathbb R (D) \oplus (\mathfrak{sl}(	q)\oplus \mathfrak{sl}(p))_C$.

\end{lemma}

\begin{remark}    Conditions are needed on the orbit to ensure the result in the lemma.  For example, for type $(2,2k+1)$ it is known that there is one generic orbit of the $GL(2k+1)\times GL(2)$ action and this orbit is open.  As such, one can see that the stabilizer of $\mathfrak{gl}(2k+1)\oplus \mathfrak{gl}(2)$ is bigger than the
 stabilizer of  $\mathfrak{sl} (2k+1)\oplus \mathfrak{sl}(2)$  extended by the (1,2) derivation given in Eqn.~\ref{eqn: 1,2 derivation}.
\end{remark}

\begin{proof}[Proof of the lemma]

We  write $\mathfrak{gl}(q)\oplus \mathfrak{gl}(p)$ as $\mathbb R \oplus \mathbb R \oplus \mathfrak{sl}(q)\oplus \mathfrak{sl}(p)$.  Each $\mathbb R$ factor acts by scaling on $C$ (cf. Eqn.~\ref{eqn: gl x gl Lie algebra action}).  Given $X\in (\mathfrak{gl}(q)\oplus \mathfrak{gl}(p))_C$, we see that it may be written as $X=X_1 + X_2$ with $X_1$ acting by scaling and $X_2\in \mathfrak{sl}(q)\oplus \mathfrak{sl}(p)$.  Thus
	$$X_2 \cdot C = r C,$$
for some $r\in\mathbb R$.  This implies 
	$$\exp(tX_2)\cdot C = e^{rt} C,$$ 
for $t\in\mathbb R$.  However, since $SL(q)\times SL(p)\cdot C$ is closed it cannot contain the origin in its boundary, thence we see that $r=0$ and  $X_2$ stabilizes $C$.  In turn, $X_1$ must stabilize $C$ and so is a multiple of the (1,2)-derivation $D$.

\end{proof}

The proof of the theorem above now follows immediately from the lemma combined with  some general theory on certain representations of $SL(q)\times SL(p)$ which we present in the following section.

\section{The action of $SL_p\times SL_q$}\label{sec: action of SL x SL}

For a detailed discussion on the following   facts and their relationship to nilpotent geometry, we refer the interested reader to \cite[Chapter 7]{Jablo:Thesis}.

From the work of Knop-Littelman \cite{KnopLittelman:DerGradErzeugenderFunktionenVonInvariantenringen}, we know that in the case of non-exceptional types $(p,q)$, the generic stabilizers of the the $SL(q)\times SL(p)$ representation above  are finite (i.e. they have trivial connected component).  From \cite{PopovBook}, we then have that generic orbits of $SL(q)\times SL(p)$ are closed.   We summarize this information below.

\begin{thm}
For non-exceptional types $(p,q)$, there is a Zariski open set in $\mathfrak{so}(q)^p$ whose elements $C$ satisfy
	\begin{enumerate}
	\item $SL(q)\times SL(p) \cdot C$ is closed, and
	\item $(\mathfrak{sl}(q)\times \mathfrak{sl}(p) )_C = \{0\}$
	\end{enumerate}
\end{thm}
Combining part (ii) of this theorem with the lemma above,  the proof of Theorem \ref{thm: generic derivations} is complete.

Determining when a particular orbit is closed is non-trivial as closedness is a global property and these orbits are non-compact.  As we will need to do this to obtain an explicit example for Theorem A, we lay out some criteria from the general theory below.

Consider the maps $m_1:\mathfrak{so}(q)^p \to symm(q)$ and $m_2:\mathfrak{so}(q)^p \to symm(p)$ defined by
	\begin{align}\label{eqn: definitions of m1 and m2}
	\begin{split}
	m_1(C) &= \sum_{i=1}^q (C_i)^2,\\
	\left(m_2\left(C\right)\right)_{ij} &= \mytrace (C_iC_j), \mbox{ for } 1\leq i,j\leq p.
	\end{split}
	\end{align}
The so-called moment map for the action of $SL(q)\times SL(p)$ on $\mathfrak{so}(q)^p$ is given by
	\begin{align}\label{eqn: moment map for sl x sl action}
	m(C) = \left( m_1(C) - \frac{\mytrace \ m_1(C)}{q} \Id , m_2(C) - \frac{\mytrace \ m_2(C)}{p} \Id \right)   .
	\end{align}
	
The elements in $m^{-1}(0)$ are called the minimal vectors of the representation and we see that
	\[ m^{-1}(0) = \left\{ C\in\mathfrak{so}(q)^p 
	\ | \ m_1(C) = \alpha \Id \mbox{ and } m_2(C) = \beta \Id,
	 \mbox{ for some } \alpha,\beta \in\mathbb R \right\}   .\]
From the general theory \cite{RichSlow}, we have the following result.

\begin{prop}
An orbit $SL(q)\times SL(p)\cdot C$ is closed if and only if $ SL(q)\times SL(p)\cdot C \cap m^{-1}(0) \not = \emptyset$.
\end{prop}

Even further, in the setting of 2-step, nilpotent geometry, these minimal vectors are precisely the metric structure constants for soliton metrics satisfying the extra geometric condition that their  Ricci tensor is invariant under the geodesic flow \cite{Eberlein:prescribedRicciTensor}.

For clarity in the following proposition, we use the following notation.
	\begin{align*}
	G &= GL(q)\times GL(p)\\
	\mathfrak g &= \mathfrak{gl}(q)\times \mathfrak{gl}(p)\\
	K &= \left\{ g\in G \ | \ g^t = g^{-1} \right\} = O(q)\times O(p)\\
	\mathfrak p &= \left\{ X\in\mathfrak g \ | \ X^t = X\right\} = symm(q)\times symm(p)
	\end{align*}

\begin{prop}\label{prop: m1 and m2 are multiples of the identity }
Take $C\in \mathfrak{so}(q)^p$.  If $m_1(C) = \alpha \Id$ and $m_2(C) = \beta \Id$ for some $\alpha,\beta \in\mathbb R$, then
	\begin{enumerate}
	\item $C$ is a minimal point of the $SL(q)\times SL(p)$ action.
	\item $SL(q)\times SL(p)\cdot C$ is closed;
	\item The stabilizer subgroup of $GL(q)\times GL(p)$
	 at $C$ is the product $G_C = K_C \ \exp(\mathfrak g_C \cap \mathfrak p)$.
	 \item When $\mathfrak g_C$ has minimal dimension, i.e. $\mathfrak g_C = \Rspan \left(\Psi\left(D\right)\right)$, one has $K_C$ is finite and $G_C = K_C \ \exp\left(\Rspan \left(\Psi\left(D\right)\right)\right)$.
	\end{enumerate}
\end{prop}
Here $D$ is the $(1,2)$-derivation defined in Equation \ref{eqn: 1,2 derivation}.  Recall, the stabilizer subalgebra and stabilizer subgroup are, respectively, isomorphic to the derivations and automorphisms that preserve the subspaces $\mathfrak v =\mathbb R^q$ and $\mathfrak z=\mathbb R^p$; see  Eqn.~\ref{eqn: isomorphism between derivations and stabilizer subalgebra} for the  definition of the isomorphism $\Psi$.

\begin{proof}
The first two items were observed by Eberlein \cite{Eberlein:prescribedRicciTensor}; we include a short proof for completeness.  

When  $m_1(C)$ and $m_2(C)$ are multiples of the identity, their traceless parts are zero and we immediately have that $m(C)=0$, cf.\ Equation \ref{eqn: moment map for sl x sl action}; that is, $C$ is a minimal point of the action $SL(q)\times SL(p)$.  The orbit $SL(q)\times SL(p)\cdot C$ is closed by the previous proposition.

Understanding the stabilizer is more subtle.  If we were considering the stabilizer of the group $SL(q)\times SL(p)$, then we would have a decomposition as above from the standard theory using minimal vectors.  In the case of the $GL(q)\times GL(p)$ action, the moment map is $m_{GL}(C) = (m_1(C),m_2(C))$ and $C$ being a minimal point for the $SL(q)\times SL(p)$ action implies
	\[ m_{GL}(C) \cdot C = \lambda C,\]
for some $\lambda \in \mathbb R$.  That is, $C$ is a critical point of the function $F = || m_{GL} ||^2$.  From the general theory \cite{Kirwan}, one knows that the stabilizer subgroup decomposes as a semi-direct product
	\[ G_C = \left(G^{\Psi(D)}\right)_C \ \exp\left( \left(\mathfrak u^{\Psi(D)}\right)_C \right), \]
where $G^{\Psi(D)}$ is the collection of elements of $G$ which commute with ${\Psi(D)}$,  while $\left(\mathfrak u^{\Psi(D)}\right)$ is the nilpotent algebra $\left\{ X\in\mathfrak g \ | \ \lim_{t\to\infty} \exp(t \Psi(D)) X \exp( t\Psi(D))^{-1} =0  \right\}$.  See \cite[Section 6]{Kirwan} for more details on the critical points of $F$, the stratification of the representation space, and the associated parabolic groups at the critical points of $F$.

As $\Psi(D)$ lies in the center of $G$, we see that $\mathfrak u^{\Psi(D)}$ is trivial.  As ${\Psi(D)}$ is symmetric, we have  $G_C = (G^{\Psi(D)})_C = K_C \  \exp(\mathfrak g_C \cap \mathfrak p)$, as desired.

The final claim follows immediately from the third as $K_C$ being compact with dimension 0  implies finite.

\end{proof}

\section{Application to nilsolitons}

It was recognized in \cite{Eberlein:prescribedRicciTensor} that if the orbit $SL(q)\times SL(p) \cdot C$ is closed, then the associated 2-step nilpotent Lie group of type $(p,q)$ admits a soliton metric.  
Applying the above facts to nilpotent geometry, we now have the following.

\begin{thm} \label{thm: generic points being soliton and small derivations}
Consider a non-exceptional type $(p,q)$.  There exists a Zariski open set of $\mathfrak{so}(q)^p$ such that the corresponding 2-step nilpotent Lie groups $N$ satisfy
	\begin{enumerate}
	\item $N$ admits a nilsoliton metric and
	\item the derivation algebra of $\mathfrak n = Lie~N$ is the minimal one described in Theorem \ref{thm: generic derivations}.
	\end{enumerate}

\end{thm}

Using any of these generic algebras, we have all the necessary conditions to carry out the recipe in Section \ref{sec: outline of proof} for building a counterexample to the local version of the conjecture by Taketomi-Tamaru.  Part of the challenge with building examples is that, in practice, one tends adds extra symmetries to make solving a geometry problem easier.  Here we are removing symmetries to have a generic object.  So, we have   the classic problem of finding the hay in the haystack.

The smallest dimension where Theorem \ref{thm: generic points being soliton and small derivations} applies is  dimension 9, and this is for algebras of type $(4,5)$.  Below,  we give an explicit example of a 9-dimensional, 2-step nilpotent Lie algebra of type $(4,5)$ which satisfies Theorem A.  Although a generic algebra of type $(4,5)$ would suffice, in practice it is very hard to know if one is holding a generic algebra.  As we explain in Remark \ref{remark: stabilizer of our explicit example}, we still do not know if our example is truly generic.

We build our example by employing the tools from Section \ref{sec: action of SL x SL}.   Consider tuple $C=(C_1,\dots , C_4)\in \mathfrak{so}(5)^4$ given by


\[
\begin{array}{cc}
C_1 = \renewcommand{\arraystretch}{2.3}
\begin{bNiceMatrix}[columns-width=1cm]
0 & -\frac{1}{2} & 0 & 0 & 0 
\\
 \frac{1}{2} & 0 & 0 & 0 & 0 
\\
 0 & 0 & 0 & -1 & 0 
\\
 0 & 0 & 1 & 0 & 0 
\\
 0 & 0 & 0 & 0 & 0 
\end{bNiceMatrix}
& 
C_2 = \renewcommand{\arraystretch}{2.3}
\begin{bNiceMatrix}[columns-width=1cm]
0 & 0 & 0 & 0 & 0 
\\
 0 & 0 & 1 & 0 & 0 
\\
 0 & -1 & 0 & 0 & 0 
\\
 0 & 0 & 0 & 0 & -\frac{1}{2} 
\\
 0 & 0 & 0 & \frac{1}{2} & 0 
\end{bNiceMatrix}
\\
C_3 = \renewcommand{\arraystretch}{2.3}
\begin{bNiceMatrix}[columns-width=1cm]
0 & 0 & 0 & 0 & \frac{3 \sqrt{2}}{4} 
\\
 0 & 0 & 0 & -\frac{\sqrt{2}}{4} & 0 
\\
 0 & 0 & 0 & 0 & 0 
\\
 0 & \frac{\sqrt{2}}{4} & 0 & 0 & 0 
\\
 -\frac{3 \sqrt{2}}{4} & 0 & 0 & 0 & 0 
\end{bNiceMatrix}
&
C_4 = \renewcommand{\arraystretch}{2.3}
\begin{bNiceMatrix}[columns-width=1cm]
0 & 0 & 0 & \frac{\sqrt{10}}{4} & 0 
\\
 0 & 0 & 0 & 0 & -\frac{\sqrt{10}}{4} 
\\
 0 & 0 & 0 & 0 & 0 
\\
 -\frac{\sqrt{10}}{4} & 0 & 0 & 0 & 0 
\\
 0 & \frac{\sqrt{10}}{4} & 0 & 0 & 0 
\end{bNiceMatrix}
\end{array}
\]

\begin{lemma}  
The point $C\in\mathfrak{so}(5)^4$, above, is a minimal vector for the $SL(5)\times SL(4)$-orbit and so $SL(5)\times SL(4)\cdot C$ is closed.  Further, the stabilizer subgroup has minimal dimension and satisfies $\left(GL(5)\times GL(4)\right)_C \supset \left(\mathbb Z_2\times \mathbb Z_2\right) \times \exp\left(\Rspan \left( \Psi\left(D\right)\right)\right)$.
\end{lemma}

\begin{remark}
Calculating the stabilizer subalgebra is a linear algebra problem.  However, finding the stabilizer subgroup is challenging.  Even with the help of a computer algebra system, we are not able to compute the full stabilizer subgroup.
\end{remark}

\begin{proof}  Recall the definitions of $m_1(C)$ and $m_2(C)$ given in Equation \ref{eqn: definitions of m1 and m2}.  Computing, we have
	\begin{align*}
	m_1(C) &= \left[\renewcommand{\arraystretch}{1.6}
	\begin{array}{rrrrr}
-2 & 0 & 0 & 0 & 0 
\\
 0 & -2 & 0 & 0 & 0 
\\
 0 & 0 & -2 & 0 & 0 
\\
 0 & 0 & 0 & -2 & 0 
\\
 0 & 0 & 0 & 0 & -2 
\end{array}\right]
\\
	m_2(C) &=\left[\renewcommand{\arraystretch}{1.6}
	\begin{array}{rrrr}
-\frac{5}{2} & 0 & 0 & 0 
\\
 0 & -\frac{5}{2} & 0 & 0 
\\
 0 & 0 & -\frac{5}{2} & 0 
\\
 0 & 0 & 0 & -\frac{5}{2} 
\end{array}\right]
	\end{align*}
By Proposition \ref{prop: m1 and m2 are multiples of the identity }, we have that $C$ is a minimal vector and the orbit $SL(5)\times SL(4)\cdot C$ is closed.

Computing the stabilizer subalgebra $\mathfrak g_C$ is a straightforward linear algebra problem that can be done by hand; however, we use Maple to double check our calculations and verify  that 
	\[ \mathfrak g_C = \Rspan (\Psi (D) ).\]
This shows that the stabilizer subalgebra and stabilizer subgroup have minimal dimension.  Lastly, we show that the stabilizer $K_C = G_C \cap \left( O(5)\times O(4) \right)$ contains a subgroup of the form $\mathbb Z_2\times \mathbb Z_2$.

Observe, for $\epsilon = \pm 1$, we have 
$\renewcommand{\arraystretch}{2}
\begin{bNiceArray}{c|c}[margin]
\epsilon \Id_5 &   \\ \hline   & \Id_4
\end{bNiceArray} \in GL(5) \times GL(4)$ lies in the kernel of the action and so lies in the stabilizer of every point.  Further, one can quickly check that 
	\[
	\renewcommand{\arraystretch}{1.4}
	\begin{bNiceMatrix}[columns-width=3mm,margin]
	\epsilon \\
	& -\epsilon\\
	&&\epsilon\\
	&&& -\epsilon\\
	&&&& \epsilon\\
	\hline
	&&&&& -1\\
	&&&&&& -1\\
	&&&&&&& 1\\
	&&&&&&&& -1
	\CodeAfter
	\tikz \draw  (1-|6) -- (10-|6); 
	\end{bNiceMatrix}
	\]
lies in stabilizer.  Clearly, these matrices all commute and we have a group isomorphic to $\mathbb Z_2\times \mathbb Z_2$.
\end{proof}

\begin{remark}\label{remark: stabilizer of our explicit example}
If one assumes that $(g,h)\in K_C$ has diagonal $h$, then $g$ must be diagonal and $(g,h)$ in subgroup $\mathbb Z_2\times\mathbb Z_2$ constructed above.  However, we do not know if $h$ must always be diagonal.  It would be interesting to know the full stabilizer subgroup of $C$.  Even further, we do not know the stabilizer in general position.  Certainly, it contains $\mathbb Z_2$ from kernel of the action.  It would be interesting to know if the stabilizer in general position must be larger.  If the generic stabilizer were only $\mathbb Z_2$, then our example would not be a generic point.  Other examples we have constructed also contain $\mathbb Z_2\times\mathbb Z_2$ in their stabilizers.
\end{remark}

\begin{example}\label{example: 9-dimensional nilpotent}
The 9-dimensional, 2-step nilpotent Lie group associated to $C$ above satisfies the conditions of Theorem A.
\end{example}

\begin{proof}
From the lemma, we know that $SL(5)\times SL(4)\cdot C$ is closed and so the 2-step nilpotent Lie algebra admits a soliton metric \cite{Eberlein:prescribedRicciTensor}.  

The derivation algebra is minimal as the stabilizer subalgebra of the $\mathfrak{sl}(5)\times \mathfrak{sl}(4)$-action has minimal dimension.  We now satisfy all the conditions of the recipe given in Section~\ref{sec: outline of proof} and we have our example.
\end{proof}

Lastly, we point out that since $q=5$ is odd, our automorphism group must be disconnected as it contains the element 
$\renewcommand{\arraystretch}{2}
\begin{bNiceArray}{c|c}[margin]
-\Id_5 &   \\ \hline   & \Id_4
\end{bNiceArray} \in GL(5) \times GL(4)$, which has negative determinant.  If $\Aut$ had only two components, then it would be a normal subgroup of $S$, cf.~discussion in Section \ref{sec: outline of proof}.  A nilpotent group with such an automorphism group would then give a full counterexample to the conjecture of Taketomi-Tamaru.  It would be interesting to know if such exists or if the automorphism group must always have more than two components - in which case it would not be a normal subgroup of $S$.

\providecommand{\bysame}{\leavevmode\hbox to3em{\hrulefill}\thinspace}
\providecommand{\MR}{\relax\ifhmode\unskip\space\fi MR }
\providecommand{\MRhref}[2]{%
  \href{http://www.ams.org/mathscinet-getitem?mr=#1}{#2}
}
\providecommand{\href}[2]{#2}


\begin{thebibliography}{GW88}

\bibitem[BL22]{Bohm-Lafuente:HomogeneousEinsteinMetricsOnEuclideanSpacesAreEinsteinSolvmanifolds}
Christoph B\"{o}hm and Ramiro~A. Lafuente, \emph{Homogeneous {E}instein metrics
  on {E}uclidean spaces are {E}instein solvmanifolds}, Geom. Topol. \textbf{26}
  (2022), no.~2, 899--936. \MR{4444271}

\bibitem[Ebe94]{Eberlein:2step}
Patrick Eberlein, \emph{Geometry of $2$-step nilpotent groups with a left
  invariant metric}, Ann. Sci. \'Ecole Norm. Sup.(4) \textbf{27} (1994), no.~5,
  611--660.

\bibitem[Ebe03]{EberleinModuli}
\bysame, \emph{The moduli space of 2-step nilpotent {L}ie algebras of type
  $(p,q)$}, Explorations in complex and Riemannian geometry, Contemporary
  Mathematics, American Mathematical Society \textbf{332} (2003), 37--72.

\bibitem[Ebe08]{Eberlein:prescribedRicciTensor}
\bysame, \emph{Riemannian 2-step nilmanifolds with prescribed {R}icci tensor},
  Geometric and probabilistic structures in dynamics, Contemp. Math., vol. 469,
  Amer. Math. Soc., Providence, RI, 2008, pp.~167--195.

\bibitem[GJ19]{GordonJablonski:EinsteinSolvmanifoldsHaveMaximalSymmetry}
Carolyn Gordon and Michael Jablonski, \emph{Einstein solvmanifolds have maximal
  symmetry}, Journal of Differential Geometry \textbf{111} (2019), no.~1,
  1--38.

\bibitem[GJ24]{GordonJablonski:RicciSolitonSolvmanifoldsHaveInfinitesimalMaximalSymmetry}
\bysame, \emph{Ricci soliton solvmanifolds have infinitesimal maximal
  symmetry}, to appear in Transactions of the American Mathematical Society
  (2024).

\bibitem[GW88]{GordonWilson:IsomGrpsOfRiemSolv}
Carolyn~S. Gordon and Edward~N. Wilson, \emph{Isometry groups of {R}iemannian
  solvmanifolds}, Trans. Amer. Math. Soc. \textbf{307} (1988), no.~1, 245--269.

\bibitem[Has14]{Hashinaga-OnTheMinimalityOfTheCorrespondingSubmanifoldsToFour-DimensionalSolvsolitons}
Takahiro Hashinaga, \emph{On the minimality of the corresponding submanifolds
  to four-dimensional solvsolitons}, Hiroshima Math. J. \textbf{44} (2014),
  no.~2, 173--191. \MR{3251821}

\bibitem[Heb98]{Heber}
Jens Heber, \emph{Noncompact homogeneous {E}instein spaces}, Invent. Math.
  \textbf{133} (1998), no.~2, 279--352.

\bibitem[HT17]{Hashinaga-Tamaru:Three-dimensionalSolvsolitonsAndTheMinimalityOfTheCorrespondingSubmanifolds}
Takahiro Hashinaga and Hiroshi Tamaru, \emph{Three-dimensional solvsolitons and
  the minimality of the corresponding submanifolds}, Internat. J. Math.
  \textbf{28} (2017), no.~6, 1750048, 31. \MR{3663797}

\bibitem[Jab08]{Jablo:Thesis}
Michael Jablonski, \emph{Real geometric invariant theory and {R}icci soliton
  metrics on two-step nilmanifolds}, Thesis (May 2008).

\bibitem[Jab11]{Jablo:ConceringExistenceOfEinstein}
\bysame, \emph{Concerning the existence of {E}instein and {R}icci soliton
  metrics on solvable lie groups}, Geometry \& Topology \textbf{15} (2011),
  no.~2, 735--764.

\bibitem[Jab19]{Jablo:MaximalSymmetryAndUnimodularSolvmanifolds}
\bysame, \emph{Maximal symmetry and unimodular solvmanifolds}, Pacific J. Math.
  \textbf{298} (2019), no.~2, 417--427. \MR{3936023}

\bibitem[Kir84]{Kirwan}
Frances~Clare Kirwan, \emph{Cohomology of quotients in symplectic and algebraic
  geometry}, Mathematical Notes 31, Princeton University Press, Princeton, New
  Jersey, 1984.

\bibitem[KL87]{KnopLittelman:DerGradErzeugenderFunktionenVonInvariantenringen}
Friedrich Knop and Peter Littelman, \emph{Der grad erzeugender funktionen von
  invariantenringen. ({G}erman) [{T}he degree of generating functions of rings
  of invariants]}, Math. Z. \textbf{196} (1987), no.~2, 211--229.

\bibitem[Lau01]{Lauret:RicciSolitonHomogeneousNilmanifolds}
Jorge Lauret, \emph{Ricci soliton homogeneous nilmanifolds}, Math. Ann.
  \textbf{319} (2001), no.~4, 715--733.

\bibitem[Nik11]{Nikolayevsky:EinsteinSolvmanifoldsandPreEinsteinDerivation}
Y.~Nikolayevsky, \emph{Einstein solvmanifolds and the pre-{E}instein
  derivation}, Trans. Amer. Math. Soc. \textbf{363} (2011), 3935--3958.

\bibitem[PV94]{PopovBook}
V.L. Popov and E.B. Vinberg, \emph{Algebraic geometry {IV}: {II}. invariant
  theory}, Springer-Verlag, Berlin Heidelberg, 1994.

\bibitem[RS90]{RichSlow}
R.W. Richardson and P.J. Slodowy, \emph{Minimum vectors for real reductive
  algebraic groups}, J. London Math. Soc. \textbf{42} (1990), 409--429.

\bibitem[Tak18]{Taketomi:OnARiemannianSubmanifoldWhoseSliceRepresentation}
Yuichiro Taketomi, \emph{On a {R}iemannian submanifold whose slice
  representation has no nonzero fixed points}, Hiroshima Math. J. \textbf{48}
  (2018), no.~1, 1--20. \MR{3771997}

\bibitem[TT18]{Taketomi-Tamaru:OnTheExistenceOfLeftInvariantRicciSolitons-AConjectureAndExamples}
Y.~Taketomi and H.~Tamaru, \emph{On the nonexistence of left-invariant {R}icci
  solitons---a conjecture and examples}, Transform. Groups \textbf{23} (2018),
  no.~1, 257--270. \MR{3763948}

\end{thebibliography}
\end{document}